\documentclass[11pt]{article}
\usepackage{amsmath,amssymb,amsfonts}

\newtheorem{theorem}{Theorem}[section]
\newtheorem{rem}[theorem]{Remark}

\newtheorem{lemma}[theorem]{Lemma}

\newtheorem{eg}{Example}
\newenvironment{example}{\begin{eg}\rm}{\end{eg}}
\newtheorem{conj}[theorem]{Conjecture}

\newtheorem{defin}[theorem]{Definition}

\input{epsf}
\author{{\normalsize{ M. Mohammad-Noori}
}\vspace{2mm} \\{\footnotesize{$^{ \textrm{a}}$\it
Department of
Mathematics, Statistics and Computer
Science, University of Tehran,}}\vspace{-2mm}\\{\footnotesize{\it   Tehran,
Iran}}
\\{\footnotesize Emails: morteza@ipm.ir, mnoori@khayam.ut.ac.ir}
}

\begin{document}
\title{Enumeration of closed random walks in the square lattice according to their areas}

 \maketitle
\begin{abstract}
We study the area distribution of closed walks of length $n$, beginning and ending at the origin.
The concept of area of a walk in the square lattice is generalized and the usefulness of the new concept
 is demonstrated through a simple argument. It is concluded that
 the number of walks of length $n$ and area $s$ equals to the coefficient of $z^s$ in the expression $(x+x^{-1}+y+y^{-1})^n$, where the calculations are performed in a special group ring $R[x,y,z]$.
  A polynomial time algorithm for calculating these values, is then concluded.
 Finally, the provided algorithm and the results of implementation
 are compared with previous works.
\end{abstract}

{\small Subject classification: 40.09, 70.04}

\section{Introduction}

The problem of finding the area distribution of random walks of a given length
is an interesting problem which has many applications, for instance in  conformations of polymers and proteins
(see \cite{minkov, Afsari} and the references there in). When $n\rightarrow \infty$ this distribution was computed first by L\'evy using Brownian paths. In \cite{belliss}, non-commutative geometry techniques are applied to the Harper equation and the asymptotic distribution of the area, enclosed by
a random walk in the square lattice, is provided. This approach is also studied in \cite{alavi}.
In \cite{mingo} the authors used a more complicated
and more informative method to derive the above asymptotic distribution.
 They have used combinatorial arguments, in which, the enumeration of
 up-down permutations and the exponential formula for cycles of permutations play fundamental roles.
 The asymptotic formula in \cite{belliss} is compared with exact results obtained by computers: For this purpose, the closed walks of lengthes $n=16$ $n=18$ and $n=20$ are enumerated according to their areas. These computations are then extended in \cite{Afsari} to $n=28$ by using algorithmic techniques and a DSP processor. However, the algorithm used there, is based on generating walks, and thus the required time grows exponentially with respect to the length of walks. Here, we present a polynomial time algorithm for the same enumeration problem.
  As examples we present the results of implementation for $n=32,64,128$ in the appendix.

  This paper is organized as follows: In Section 2, we generalize the concept of area for a desired (not necessarily closed) walk in the square lattice. We see that the area of composition of two walks satisfies a nice relation and this approves the usefulness of the generalized concept. In Section 3, we show that previous generalization naturally lead to algebraic tools, namely a generating function $\mathfrak{w}^n=(x+x^{-1}+y+y^{-1})^n$ and a group ring $R[x,y,z]$ which together simplify the process of computation of the area distribution. This is explained in Theorem \ref{enumerative2} and the details following it. The algorithm is then discussed in Section 4 and some selected
  results obtained from the implementation are given in tables in the appendix.

\section{Area of walks in the square lattice}
It is well-known that the algebraic area corresponding to a
polygon ${\bf P}=A_0$\\ $A_1 \ldots A_{n-1}$ with coordinates
$A_i=(x_i,y_i), 0\leq i\leq n-1$ and $(x_n,y_n)=(x_0,y_0)$, is as
follows
\begin{equation}\label{AreaPolyg}S({\bf P})=\frac{1}{2}\sum_{i=0}^{n-1}(x_i y_{i+1} - x_{i+1} y_i).\end{equation}

If moreover, the sides of the polygon are parallel to the axes $x$ or $y$ i.e. if the identity
\begin{equation}\label{conditns}
(x_i-x_{i+1})(y_i-y_{i+1})=0,
\end{equation}

holds for $0\leq
i \leq n-1$, then
\begin{equation}\label{Area1} S({\bf P})=\sum_{i=0}^{n-1} x_i(y_{i+1}-y_i)\ .\end{equation}

By setting $\epsilon_i=y_{i+1}-y_i$ and $\delta_j=x_{j+1}-x_j$, it is obtained that

\begin{equation}\label{Area11} S({\bf P})= \sum_{0\leq j< i \leq n-1} \epsilon_i \delta_j\ .\end{equation}

 Now, consider a more general case where we have just a sequence of
 vertices ${\bf Q}=A_0 A_1 \ldots A_n$ in which the first
 and the last points are not necessarily the same, but the coordinates satisfy condition
(\ref{conditns}) for $0\leq i\leq n-1$.

\noindent In this case (\ref{Area11}) can be defined as the area of ${\bf Q}$. To rewrite this formula in terms of
coordinates $x_i$ and $y_i$, we do as follows:
\begin{align*}
S({\bf Q})&= \sum_{i=0}^{n-1}\epsilon_i \sum_{j=0}^{i-1}\delta_j\\
&= \sum_{i=0}^{n-1}\epsilon_i(x_i-x_0)
\end{align*}
and since $\sum_{i=0}^{n-1}\epsilon_i=y_n-y_0$, the area is obtained as
 \begin{equation}\label{AreaOpen} S({\bf Q})=\sum_{i=0}^{n-1}{x_i(y_{i+1}-y_i) + x_0 (y_0 - y_n)}.\,\end{equation}

This equation has a simple geometric interpretation: The area
of ${\bf Q}$ equals the area of the following polygon
$${\bf Q}_c=[(x_0,y_0),(x_1,y_1),\ldots,(x_n,y_n),(x_0,y_n),(x_0,y_0)].$$

Now, consider a walk beginning at $(x_0,y_0)$ in the square lattice. Each step is just moving
one unit to right, left, up or down. We may denote these moves
simply by $r$, $\bar{r}$, $u$, and $\bar{u}$. Thus a walk of
length $n$ may be demonstrated by its first point $(x_0,y_0)$ and a word $w=w_0\ldots w_{n-1}$ of
length $n$ over the alphabet $\mathcal{A}=\{r,\bar{r},u,\bar{u}\}$.
This walk can alternatively be demonstrated by corresponding sequence of vertices
 $$W=[(x_0,y_0),(x_1,y_1),\ldots,(x_n,y_n)],$$
 which not only satisfies condition \ref{conditns}, but also the following equalities
 for $0 \leq i \leq n-1$ :
\begin{equation} \label{HvyCond}\{|x_{i+1}-x_i|,|y_{i+1}-y_i|\}=\{0,1\}.\end{equation}

Thus the equation (\ref{AreaOpen}), can naturally be used for the area of a walk in ${\mathbb{Z}}^2$. The area of such a walk is obviously independent from the beginning point $(x_0,y_0)$ and depends only on the word $w$.

 Two walks can simply be composed as follows: Let
$$W=[(x_0,y_0),(x_1,y_1)\ldots,(x_{n},y_{n})],$$
$$W'=[(x'_0,y'_0),(x'_1,y'_1),\ldots,(x'_{m},y'_{m})].$$
Then we define the composed walk $WW'$ as
$$ WW'=[(x_0,y_0),\ldots,(x_{n},y_{n}),(x_{n+1},y_{n+1}),\ldots,(x_{n+m},y_{n+m})],$$
where $x_{n+i}=x_{n}+x'_i-x'_0$ and $y_{n+i}=y_{n}+y'_i-y'_0$ for
$1 \leq i \leq m$. This composition corresponds to concatenation of words $w$ and $w'$. It is easy to prove that
$$S(WW')=S(W)+S(W')+(x_{n}-x_0)(y'_{m}-y'_0).$$

 Let $(x_0,y_0)=(x'_0,y'_0)=(0,0),
(x_{n},y_{n})=(i,j),(x_{m},y_{m})=(i',j')$. Then
\begin{equation} \label{SWWprim}S(WW')=S(W)+S(W')+ij'.\end{equation}

The geometric interpretation of this fact is demonstrated in Figure 1.

\vspace*{.1cm}
{\small $$\epsfbox{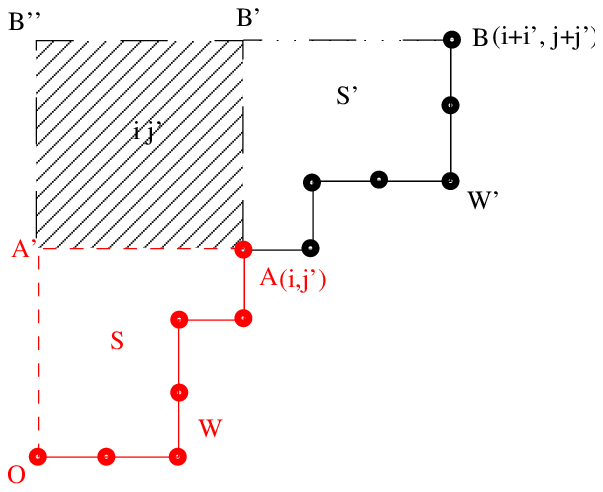}$$}
\begin{center}
{\rm Figure 1}
\end{center}
\vspace*{.1cm}

\section{A generating function, a new multiplication and algorithmic results}

The following enumerative result about the number of
walks in the square lattice is well-known and easy to prove :

\begin{lemma}\label{enumerative1} Let $\mathfrak{u}=\alpha+\alpha ^{-1}+\beta+\beta^{-1}$ and
denote the number of walks in the square lattice which begin at
the origin and end in $(p,q)$ by $a_n(p,q)$. Then

\begin {itemize}
\item[(i)] {The number $a_n(p,q)$ equals the coefficient of $\alpha^p
\beta^q$ in $\mathfrak{u}^n$.} \item[(ii)]{We have $a_n(p,q)={n \choose
{\frac{n+p+q}{2}}}{n \choose {\frac{n-p+q}{2}}}$.}
\end{itemize}\end{lemma}

{\bf Notation. } In this paper, we use the notation $[x^iy^jz^s]F(x,y,z)$, to show the coefficient of $x^iy^jz^s$ in $F(x,y,z)$. Thus Lemma \ref{enumerative1}(i), states that
  $$a_n(p,q)=[\alpha^p \beta^q]\mathfrak{u}^n. $$

We remark that part (ii) of the above lemma can be proved either as an immediate consequence
of part (i) or independently through a combinatorial argument by projecting the walks on axes $y=x$ and $y=-x$ (This technique is well-known, see for instance Proposition 2.3 of \cite{mingo} and page 451 of \cite{Lot3}).

Now, let $w_n(i,j,s)$ be the number of walks
beginning at the origin and ending in $(i,j)$ and having algebraic
area $s$. Then by using identity (\ref{SWWprim}) we have the following enumerative identity:
\begin{equation}\label{recursive} w_{n+m}(i,j,s)=\sum_{I_{ijs}} w_n(i_1,j_1,s_1)w_m(i_2,j_2,s_2),\,  \end{equation}
where $I_{ijs}$ consists of the set of pairs of integer triples
 $(i_1,j_1,s_1)$ and $(i_2,j_2,s_2)$ which satisfy the following set of equations:
\begin{equation} \label{sumLim} i_1+i_2=i,\ j_1+j_2=j,\ s_1+s_2+i_1j_2=s. \end{equation}

Equations (\ref{recursive}) and (\ref{sumLim}) lead us to study
the multiplication of two monomials $X= x^i y^j z^s$ and
$X'=x^{i'} y^{j'} z^{s'}$ defined as
\begin{equation} \label{gpDef} x^i y^j z^s.x^{i'} y^{j'} z^{s'}= x^{i+i'}
y^{j+j'} z^{s+s'+ij'}
 \end{equation} It is easy to check that the set of all
monomials $x^i y^j z^s$ with this multiplication construct a
non-commutative group. Note that the element $x^i y^j z^s$ is just
a representation of an element of a group (which may also be
represented just by the triple $(i,j,s)$) and in general does not
equal $x^i.y^j.z^s$ (for instance $x.y.z=xyz^2$). This leads to
construct a group ring with elements $\sum a(i,j,s) x^i y^j z^s$
with finitely many nonzero coefficients which come from a ring, say
$\mathbb{R}$. Now considering the group ring $R[x,y,z]$, the
following theorem is a generalization of Lemma \ref{enumerative1}:

\begin{theorem}\label{enumerative2} Let $\mathfrak{w}=x+x^{-1}+y+y^{-1}$ and denote the number of walks in the
square lattice which begin at the origin and end in $(p,q)$ and
has area $s$ by $w_n(p,q,s)$. Then
 {The number $w_n(p,q,s)$ equals the coefficient of $x^p y^q z^s$ in $\mathfrak{w}^n$, where computations
  are performed in $R[x,y,z]$: $$w_n(p,q,s)=[x^p y^q z^s]\mathfrak{w}^n $$.}
\end{theorem}

Obviously the mapping $\tau:R[x,y,z]\rightarrow R[\alpha,\beta]$ defined by $\tau(x^i y^j z^s)=\alpha^i \beta^j$ constructs a ring homomorphism.
 Note that calculating an expression of the form $(x+x^{-1}+y+y^{-1})^n$ is an idea, which is not only common
 between theorem \ref{enumerative2} and Lemma \ref{enumerative1}, but also it is used in
 \cite{belliss}. The important point here, is that in our approach, the function $\mathfrak{w}$ and the group ring $R[x,y,z]$, both appear as a conclusion of a purely combinatorial discussion by a proper extension of the concept of area for open walks.

\begin{example}\label{W2}
It is easily checked that {\scriptsize
\begin{align}
\mathfrak{w}^2 = &\ (x+x^{-1}+y+y^{-1}).(x+x^{-1}+y+y^{-1}) \nonumber\\
= & x^2+1+xy+xy^{-1}+1+x^{-2}+x^{-1}y+x^{-1}y^{-1} \nonumber\\
& + xyz+x^{-1}yz^{-1}+y^2+1+xy^{-1}z^{-1}+x^{-1}y^{-1}z+1+y^{-2} \nonumber\\
=&x^2+x^{-2}+y^2+y^{-2}+4+xy+xy^{-1}+x^{-1}y+x^{-1}y^{-1} \nonumber\\
&+xyz+x^{-1}yz^{-1}+xy^{-1}z^{-1}+x^{-1}y^{-1}z. \nonumber
\end{align}
} \noindent The terms of the form $axyz^i$ are $xy$ and $xyz$.
This means that there are just two walks of length $2$ which begin from the origin and end in $(1,1)$:
 one walk with area $0$ and another one with area $1$.
\end{example}

\begin{example}\label{W4}
Again consider the previous example. If we want to answer the same
question about walks of length $4$, it is enough to calculate the
terms $axyz^i$ in $\mathfrak{w}^4$. But using $\mathfrak{w}^4=(\mathfrak{w}^2)^2$ and applying the
above result, these terms are as follows:
\begin{align}
 &\ x^2.x^{-1}y+x^2.x^{-1}yz^{-1}+ y^2.xy^{-1}+y^2.xy^{-1}z^{-1}+4xy+4xyz+\nonumber\\
 &\ 4xy+xy^{-1}.y^2+x^{-1}y.x^2+4xyz+x^{-1}yz^{-1}.x^2+xy^{-1}z^{-1}.y^2\nonumber\\
=&\ xy(2z^2+10z+10+2z^{-1}). \nonumber
\end{align}

\end{example}

\section{An algorithm to enumerate walks}

Similar to examples \ref{W2} and \ref{W4} one can compute the expression $\mathfrak{w}^n$ for given values of $n$ to obtain values $w_n(i,j,s)$. Of course this can be done for any positive integer $n$ (not only powers of $2$,) by calculating expressions of form $\mathfrak{w}^{n_1}.\mathfrak{w}^{n_2}$ at most $\lg (n)$ times.
To analyze the provided algorithm, first note that
if $|i|+|j|>n$ or if $n+i+j$ is odd or if $|s|>\frac{n^2}{4}$, then $w_n(i,j,s)=0$ (More generally, if $|s|>\frac{(n+|i|+|j|)^2}{16}$ then $w_n(i,j,s)=0$). Thus the expression $\mathfrak{w}^n$ has $O(n^4)$ nonzero terms and computation of the expression $\mathfrak{w}^{n_1}.\mathfrak{w}^{n_2}$ needs totally $O({n_1}^4 {n_2}^4)$ multiplications. Thus for calculating $\mathfrak{w}^{n}$, the number of required integer multiplications is $O(n^8 \lg(n))$. However, as $n$ grows larger, the coefficients $w_n(i,j,s)$ grow exponentially and should be considered as ``large integers" (instead of integers) and the integer multiplication should not be considered as a unitary operation.
 (For computation with large integers, see for instance \cite{algor}).
  This problem is resolved by using modular arithmetic as follows: Choose $k$ and distinct prime integers $p_1,\cdots,p_k$ such that $w_n(i,j,s)<p_1\cdots p_k$ (It is enough to select these numbers such that $p_1\cdots p_k>4^n$). For any $i$, $1\leq i\leq k$, calculate the coefficients of $\mathfrak{w}^{n}$ $\mod p_i$. Finally for each $s$ with $s\leq n^2/16$ a sequence $t_1,\cdots,t_k$ with $w_n(0,0,s)\equiv t_i\, (\mod p_i)$ is obtained. Thus the values $w_n(0,0,s)$ can easily be reconstructed using Chinese reminder theorem. Since $k=O(\lg (n))$, the complexity of our algorithm in this case (i.e. when $4^n$ is a large integer), is computed as $O(n^8 \lg^2(n))$.

We have implemented this for $n=8, 16, 32, 64$; moreover we have obtained the terms $w_{128}(0,0,s)$ in the expression $\mathfrak{w}^{128}$ (As mentioned before, since the coefficients are large for $n=64, 128$, we have used some modular arithmetic to simplify the calculations). For any even integer $n$, $w_n(0,0,s)$ is a unimodal sequence which is symmetric with respect to the axis $s=0$ and takes its maximum at $s=0$.  The results of computations are briefly demonstrated in tables 1.1, 1.2, 2.1, 2.2 (Due to the symmetry, negative values of $s$ are omitted from these tables). Histograms of the number of closed walks with corresponding areas are demonstrated in table 1.1 for $n=16,32,64$ and in table 1.2 for $n=128$. We have used sampling of areas to estimate the whole distribution of walks with respect to their area in tables 2.1 for $n=32, 64$ and in table 2.2 for $n=128$.
We remark that it is possible to improve the complexity of this algorithm to $O(n^6\lg^3(n))$ by some modifications.
However, we do not need to implement this modified version for $n\leq 128$.


{\bf Acknowledgement.} The author would like to thank J. Bellissard for sending him the reference \cite{mingo}. Also he would like to thank J. Zahiri, N. Sediq and Sh. Mohammadi for studying a draft of this manuscript.

\newpage

\noindent{\large Table $1.1$}. Histogram of the number of closed walks for $n= 16, 32, 64$.
{\scriptsize
\begin{center}
\begin{tabular}{|r|r|r| r|}
\hline
 {Area} $s$ &\ $n=16$ &$n=32$& $n=64$\\[0.2cm]
\hline
 $0$  & $33820044$ & $3581690,9974343308$  & $165545,3003285874,5794673311,4483378060$\\
 $1$  & $28133728$  & $3444028,9607452416$ & $163951,1042472083,1102818268,2298389120$ \\
 $2$  & $18569808$ &  $3073077,4567275040$ & $159289,6955895232,9652405885,9706370944$ \\
 $3$  & $10127744$  &  $2565083,0257099008$ & $151905,3681689474,6366542956,7889122560$\\
 $4$  & $5015108$  &  $2025070,2695492528$ & $142312,9428052661,5202704645,6952968352$\\
 $5$  & $2289760$ & $1529170,0875844800$ & $131125,1853389455,8748104447,7490274432$ \\
 $6$  & $1036368$ & $1116254,1356438464$ & $118976,8131598077,3546151378,9923455744$\\
$7$  & $435040$ & $794324,3235665408$ &$106459,4734512548,5432356922,6277637888$\\
$8$  & $184104$ & $554823,8812436036$ &$94076,4049114445,2830982857,7473179632$\\
$9$  & $73056$ & $382197,6073766784$ &$82219,0388885868,0062962560,1052407296$\\
$10$  & $28064$ & $260613,1312522976$ &$71162,7986902081,4394878378,0740454720$\\
$11$  & $10336$ & $176316,8848622336$ &$61076,6841799750,4289859727,5624957312$\\
$12$  & $3760$ & $118576,4173049648$ &$52040,7460865277,2477622381,7129090240$\\
$13$  & $1088$ & $79335,0438261504$ &$44066,4566930747,9914521591,0107701504$\\
$14$  & $352$ & $52859,9386326560$ &$37116,4869321275,0343143165,4198093312$\\
$15$  & $96$ & $35083,2035517248$ &$31121,8811662498,3702411709,7801894016$\\
$16$  & $16$ & $23202,5420494728$ &$25995,8037652707,6596386739,4450277316$\\
$17$  & $0$ & $15290,8309279936$ &$21643,8097824910,3762445371,2586987264$\\
$18$  & $0$ & $10044,4272588768$ &$17971,0315275440,8547703288,0937254400$\\
$19$  & $0$ & $6574,8927845440$ &$14886,8517633991,1765417370,3945559808$\\
$20$  & $0$ & $4289,6632260736$ &$12307,6599368221,5520959165,4379143840$\\
$21$  & $0$ & $2788,4857169280$ &$10158,2235010504,6499042413,6544733440$\\
$22$  & $0$ & $1806,4537994848$ &$8372,1122544237,1269186309,3162306688$\\
$23$  & $0$ & $1165,8943874752$ &$6891,5113696382,9340030299,8512645376$\\
$24$  & $0$ & $749,7160572048$ &$5666,6701469563,7412442711,0537836464$\\
$25$  & $0$ & $480,1211445312$ &$4655,1587087675,5085653450,8420517888$\\
$26$  & $0$ & $306,2945599680$ &$3821,0488492849,7035370375,8706331264$\\
$27$  & $0$ & $194,5755843520$ &$3134,0926012667,5420527449,4192796416$\\
$28$  & $0$ & $123,0912937696$ &$2568,9430707074,7787699462,7648267040$\\
$29$  & $0$ & $77,5044394624$ &$2104,4412234929,2181969190,3419356416$\\
$30$  & $0$ & $48,5883898144$ &$1722,9792841968,7831969876,7886027968$\\
$31$  & $0$ & $30,3067180160$ &$1409,9425700569,3405222424,1279150592$\\
$32$  & $0$ & $18,8158770672$ &$1153,2267976836,0220074562,9419251848$\\
$33$  & $0$ & $11,6190755520$ &$942,8248707710,9783233726,8937696128$\\
$34$  & $0$ & $7,1372120768$ &$770,4760703796,4923852241,8443915072$\\
$35$  & $0$ & $4,3588560640$ &$629,3700506407,3095007510,8661437568$\\
$36$  & $0$ & $2,6468754368$ &$513,8983895305,2165176788,5718002464$\\
$37$  & $0$ & $1,5971326400$ &$419,4468635943,1128763468,3935688320$\\
$38$  & $0$ & $9580227072$ &$342,2223782721,8792023636,1272721856$\\
$39$  & $0$ & $5704976448$ &$279,1091340556,7526383551,8177452288$\\
$40$  & $0$ & $3374362720$ &$227,5493796993,4261990343,2651390528$\\
$41$  & $0$ & $1979897600$ &$185,4447127237,3591685480,3318881920$\\
$42$  & $0$ & $1153531776$ &$151,0745274528,2500466045,4542766592$\\
$43$  & $0$ & $665930496$ &$123,0287071199,9349261843,0645507456$\\
$44$  & $0$ & $381403552$ &$100,1521446247,9515175280,4221836704$\\
$45$  & $0$ & $216272192$ &$81,4990481851,9987650971,6280904448$\\
$46$  & $0$ & $121397120$ &$66,2953446385,4157940925,1272815296$\\
$47$  & $0$ & $67391168$ &$53,9077623054,4736632512,9024635136$\\
$48$  & $0$ & $37007392$ &$43,8184278757,0688296731,2469700240$\\
$49$  & $0$ & $20046912$ &$35,6040025857,9755520652,8140275072$\\
$50$  & $0$ & $10730048$ &$28,9185589829,9232653138,8592603584$\\

\hline

\end{tabular}

\end{center}
}
\newpage
\noindent{\large Table $1.2$}. Histogram of the number of closed walks for $n= 128$.
{\scriptsize
\begin{center}
\begin{tabular}{|r|r|}
\hline
 {Area } $s$ &\ $n=128$ \\[0.2cm]
\hline
 $0$  & $1410,7033892003,4556275957,3855536443,1713372583,8745556276,5835946782,8699656588$\\
 $1$  & $1407,3024540489,7561178225,1492421016,5903644384,0838340709,2111460482,4937387520$ \\
 $2$  & $1397,1649733707,8547470736,0499774263,8920388680,8359080398,3163318225,4162311040$ \\
 $3$  & $1380,4842678848,7302789379,1645965215,5579137301,3499661185,3647061656,3763726848$\\
 $4$  & $1357,5738773060,9946518605,0585703091,4103655038,8810690309,8680556206,4476333696$\\
 $5$  & $1328,8551923721,7677597822,7209730464,4139815783,0161065810,1704581134,1916138496$ \\
 $6$  & $1294,8413214797,3192840866,6424244840,7914453824,1990913897,2602303118,5542646784$\\
$7$  & $1256,1182279806,0973554574,6222544230,8237684539,6052335327,6433750005,5423763456$\\
$8$  & $1213,3242855140,0028030693,5376370151,6590699200,7160036130,3615131804,7547629472$\\
$9$  & $1167,1294078548,3705936097,4700188169,9134665437,9964875422,2126644178,1188844032$\\
$10$  & $1118,2148257931,4876464340,5178198327,6272565882,5003428615,7423056334,4507678592$\\
$11$  & $1067,2544253010,0491533514,7894490482,6439471597,5277883548,2780506254,9148335104$\\
$12$  & $1014,8983533104,4634867110,1206819214,1193382594,8253734997,0143386378,9088034176$\\
$13$  & $961,7593661773,7189043046,3337315981,8652756444,2697345733,0322031818,4291592704$\\
$14$  & $908,4021657996,3631929848,5392789504,0553785827,1399264505,9083287686,8407903488$\\
$15$  & $855,3357592722,7904570583,4071833016,3005561629,6188755045,2681666035,3268214784$\\
$16$  & $803,0087037290,4609496133,5353158915,5592171685,1608553021,1235191565,6504449136$\\
$17$  & $751,8069661057,2813072696,6629583064,1148708036,7436477079,4895846268,6542646784$\\
$18$  & $702,0540396862,5146881526,9860430186,5491459096,2468812491,6045345873,6116161536$\\
$19$  & $654,0129126309,4424638675,1129960649,6515635171,2836056527,8268318381,3914975488$\\
$20$  & $607,8894723270,7517045571,1069505624,9654909825,9691568380,4579410665,3886321856$\\
$21$  & $563,8369457848,1328923268,0251406724,5624510978,4465364091,8338626540,2448514048$\\
$22$  & $521,9610124403,1971690875,7285266312,1361048967,8464554560,1550600804,2947782528$\\
$23$  & $482,3252741185,5977063374,9204570737,0096681055,7983960514,9055938353,6273292544$\\
$24$  & $444,9568211123,2677528406,7913521950,2168637496,3514310616,1020629339,0696154560$\\
$25$  & $409,8516882597,8820945433,0086445538,6875648834,3075383182,7515093856,2496782336$\\
$26$  & $376,9800468596,7878523294,6468408515,6838053818,2117951574,8956155126,1898202880$\\
$27$  & $346,2910248525,7910212051,9918713267,6285219824,3417736320,7150365881,6662952192$\\
$28$  & $317,7170876122,6649370409,8769310919,3372955110,9885406167,7375795751,0850541184$\\
$29$  & $291,1779444806,8780735863,2682558559,3943314167,7299990195,1784320302,7017344000$\\
$30$  & $266,5839720083,3217144183,2330534321,5631455953,9593549729,2939952449,8739436928$\\
$31$  & $243,8391642983,0152229604,1800095915,7149399101,9378288422,7765033861,2981857280$\\
$32$  & $222,8436346973,8653725051,1841404543,0931745709,4233102150,3754802180,2932444612$\\
$33$  & $203,4957022281,8650274258,8409511283,0795000912,6351287804,8241811701,0044212480$\\
$34$  & $185,6936015120,1442388100,3756227758,1997286304,1641800578,0316304872,5931503104$\\
$35$  & $169,3368573323,1966937758,8587669558,1677928523,5336728042,7009350227,8754297344$\\
$36$  & $154,3273651764,1673709820,8254120503,0091312669,4620532769,1964334067,6582675712$\\
$37$  & $140,5702177017,1811870033,4744950358,1126325851,8097400944,1897482027,1138028544$\\
$38$  & $127,9743146199,7216121242,6590547944,5795679103,3517158052,5703219655,5346880768$\\
$39$  & $116,4527903905,2565093234,5991490752,6931680721,8392377764,6639015944,4253611008$\\
$40$  & $105,9232906770,0813934257,9389349588,9826312713,6299361111,8150382831,9241859104$\\
$41$  & $96,3081249848,4003726572,9985651294,9767291595,5105092725,2887050944,2661204480$\\
$42$  & $87,5343194257,0689809870,9084153396,1633523208,5765141173,2963363894,5553767168$\\
$43$  & $79,5335902615,1001088786,2229069228,4944928161,2471305895,5160229025,3471878400$\\
$44$  & $72,2422558339,9992080512,9735364612,5415047397,5072102242,8523062539,9947401984$\\
$45$  & $65,6011017251,4946098123,3635319649,9707716860,2623541436,0238379304,1379316736$\\
$46$  & $59,5552115318,2629696048,9600650535,3649769161,7005888983,0456617439,1624406784$\\
$47$  & $54,0537734745,4545929513,5620500429,2190688647,1817158178,7000139728,0834058496$\\
$48$  & $49,0498711810,5005814800,2302396044,8526299819,6640254673,3055286303,8011851696$\\
$49$  & $44,5002653713,1010783245,2854022296,7683798570,5477576649,4078144220,8843469312$\\
$50$  & $40,3651717975,4268971940,8275974742,9828479525,6161271917,9464626860,0827549440$\\

\hline
\end{tabular}
\end{center}

 \vspace*{0.4cm} }

{\large Table $2.1$}. Values of $w_{32}(0,0,2k)$ and $w_{64}(0,0,8k)$
{\scriptsize
\begin{center}
\begin{tabular}{|r|r|r| r|}
\hline
 {Area} $s$ &\ $n=32$ &Area $s$ &\ $n=64$\\[0.2cm]
\hline
 $0$  & $35816909974343308$ & $0$  & $165545,3003285874,5794673311,4483378060$\\
 $2$  & $30730774567275040$  & $8$ &  $94076,4049114445,2830982857,7473179632$ \\
 $4$  & $20250702695492528$ &  $16$ & $25995,8037652707,6596386739,4450277316$ \\
 $6$  & $11162541356438464$  &  $24$ & $5666,6701469563,7412442711,0537836464$\\
 $8$  & $5548238812436036$  &  $32$ & $1153,2267976836,0220074562,9419251848$\\
 $10$  & $2606131312522976$ & $40$ & $227,5493796993,4261990343,2651390528$ \\
 $12$  & $1185764173049648$ & $48$ & $43,8184278757,0688296731,2469700240$\\
$14$  & $528599386326560$ & $56$ & $8,2365472231,4016874589,9460134624$\\
$16$  & $232025420494728$ & $64$ &$1,5096688485,6768973162,6864590576$\\
$18$  & $100444272588768$ & $72$ &$2694437402,9152872362,8592916864$\\
$20$  & $42896632260736$ & $80$ &$467552007,4041838163,0604620576$\\
$22$  & $18064537994848$ & $88$ &$78744016,9083977884,8133019296$\\
$24$  & $7497160572048$ & $96$ & $12846824,7250409523,5350912480$\\
$26$  & $3062945599680$ & $104$ &$2025952,5183890308,8567109088$\\
$28$  & $1230912937696$ & $112$ &$308080,5223088698,5702239872$\\
$30$  & $485883898144$ & $120$ &$45051,1484837085,3522061280$\\
$32$  & $188158770672$ & $128$ &$6315,2306732457,4256973920$\\
$34$  & $71372120768$ & $136$ &$845,5556839591,3528941504$\\
$36$  & $26468754368$ & $144$ &$107,6805280447,7651259776$\\
$38$  & $9580227072$ & $152$ &$12,9785607608,7995391104$\\
$40$  & $3374362720$ & $160$ &$1,4718233876,4286499776$\\
$42$  & $1153531776$ & $168$ &$1559345348,1572357568$\\
$44$  & $381403552$ & $176$ &$153004502,2583865088$\\
$46$  & $121397120$ & $184$ &$13753291,0981167232$\\
$48$  & $37007392$ & $192$ & $1116802,6179713536$\\
$50$  & $10730048$ & $200$ & $80423,8839635904$\\
$52$  & $2932896$ & $208$ &$5007,3152157248$\\
$54$  & $743168$ & $216$ & $259,8435002240$\\
$56$  & $172224$ & $224$ & $10,6187399552$\\
$58$  & $35392$ & $232$ &$3102345664$\\
$60$  & $5984$ & $240$ &$53445952$\\
$62$  & $704$ & $248$ &$337280$\\
$64$  & $32$ & $256$ &$64$\\
\hline
\end{tabular}
\end{center}

 \vspace*{0.4cm} }
\newpage

\noindent{\large Table $2.2$}. Values of $w_{128}(0,0,32k)$
{\scriptsize
\begin{center}
\begin{tabular}{|r|r|}
\hline
 {Area} $s$ &\ $n=128$ \\[0.2cm]
\hline
 $0$  & $1410,7033892003,4556275957,3855536443,1713372583,8745556276,5835946782,8699656588$ \\
 $32$  & $222,8436346973,8653725051,1841404543,0931745709,4233102150,3754802180,2932444612$ \\
 $64$  & $10,1648073494,4923419172,3106847671,7720107384,7236429835,9359307733,6892311304$ \\
 $96$  & $4138993608,0697485547,8360086674,4923420517,9089245541,4042051861,9613801232$ \\
 $128$  & $161039050,5304195851,5419275033,5083133903,1254858120,1415241780,7554992752$ \\
 $160$  & $5990655,6605275424,0970303298,7170349389,5743855708,8743098925,7949721120$ \\
 $192$  & $212607,0814985186,5951934575,3598563447,3878472768,9106908764,3180434272$ \\
$224$  & $7180,6854685731,4031517786,4997582279,9706234195,6090873096,6015299584$ \\
$256$  & $230,1750942525,0140859435,7574600133,2584273010,4614639894,3043516640$ \\
$288$  & $6,9816042224,1138277250,0582452617,1692631908,7011142507,9215100160$ \\
$320$  & $1997190758,6553543306,5729674881,2193385181,4767886062,4658510144$ \\
$352$  & $53685272,9803391380,6122089350,6109205977,4658578703,0650139264$ \\
$384$  & $1350455,5535411518,0818176830,9757666811,6769796180,3565932544$ \\
$416$  & $31644,2527945050,2410030323,1585016453,7594319542,2645555776$ \\
$448$  & $687,1315164172,5008635859,1936569312,0498550867,2721950208$ \\
$480$  & $13,7450595212,7015217478,1091307692,8909430481,5008168896$ \\
$512$  & $2515774906,7368382551,6728163527,5448943450,0831010624$ \\
$544$  & $41803611,2117105753,9735632910,5539623620,6960007680$ \\
$576$  & $624889,6425097710,2464041926,3546679600,4985488768$ \\
$608$  & $8312,9826788959,9447485149,7366840822,7746173056$ \\
$640$  & $97,1577615295,3509138500,6659985912,5226822016$ \\
$672$  & $9821475660,0178800451,7922679241,8657940480$ \\
$704$  & $84231795,4751454140,9292439100,2628160000$ \\
$736$  & $598159,1733186708,6431903566,6217660160$\\
$768$  & $3408,4164367454,1023022248,2423726336$ \\
$800$  & $14,9447433526,3720269236,4107405952$ \\
$832$  & $475760748,0371846789,7834797824$ \\
$864$  & $1010220,8082355332,7727527680$ \\
$896$  & $1253,0843615061,8578411264$ \\
$928$  & $7227374040,4149925760$ \\
$960$  & $1233125,5335552896$ \\
$992$  & $19,8894005504$ \\
$1024$  & $128$ \\

\hline
\end{tabular}
\end{center}

 \vspace*{0.4cm} }

\end{document}